\documentclass[10pt,a4paper,times]{article}
\usepackage[latin1,applemac]{inputenc}
\usepackage{amsmath, amssymb, amsfonts, amsthm, amscd, epsfig, multicol, marvosym, float,mathrsfs,enumerate,appendix,float}
\usepackage[left=2cm,top=3cm,bottom=3cm,right=2cm]{geometry}
\usepackage[scaled=.95]{helvet} 
\usepackage{courier} 

\numberwithin{equation}{section}

\begin{document}
\title{A hypergeometric treatment to explain the nonlinear true behavior of redundant constraints on a straight elastic rod}
\author{Giovanni Mingari Scarpello \footnote{giovannimingari@yahoo.it}
\and Daniele Ritelli \footnote{Dipartimento scienze statistiche, viale
Filopanti, 5 40127 Bologna Italy, daniele.ritelli@unibo.it}}
\date{}
\maketitle

\begin{abstract}
In theory and practice of elastic straight rods, the statically indeterminate reactions acted by perfect constraints are commonly believed not to depend on the flexural stiffness $EJ$.
We solve exactly two elastica problems in order to obtain hypergeometrically (helped by Lagrange, Lauricella, Appell), the true displacements upon which the forces method is founded.
As a consequence, the above reactions are found to depend on stiffness: the presumptive independence credited as general, is far from being always true, but, quite the contrary, is valid only within a first-order approximation.

\vspace{0.5cm}


\end{abstract}

\section{Introduction}

In 1691, Jakob Bernoulli proposed to find out the deformed
centerline (planar "elastica") of a thin, homogeneous, straight and
flexible rod under a force applied at its end. He established, \textit{Curvatura
laminae elasticae} (1694), the differential equation of the centerline, say $y=y(x)$, of the rod \textit{ab
appenso pondere curvata} to be: 
\[
\mathrm{d}y=\frac{x^{2}\mathrm{d}x}{\sqrt{1-x^{4}}} 
\]
His nephew, Daniel Bernoulli, computed the potential energy
stored in a bent rod and, in 1742, wrote to L. Euler that the rod should
attain such a shape as to minimize the functional of squared curvature.
Accordingly, Euler dealt elastica as an isoperimetric problem of Calculus
of Variations, \textit{De curvis elasticis} (1744), arriving at elastica's
differential equation and identifying nine shapes of the curve, all equilibrium states of minimal energy. 

Anyway, C. Maclaurin had realized the elastica's connection to elliptic integrals, see 
\textit{A treatise on fluxions} (1742), \S\ 927 entitled: ``\textit{The
construction of the elastic curve, and of other figures, by the
rectification of the conic sections''.} In 1757, Euler authored a paper, 
\textit{Sur la force des colonnes}, concerning the buckling of columns
again, where the critical load is approached through a simplified expression
to the elastica's curvature.

Some analytical solutions to elastic planar curves through elliptic integrals
of the first and second kind, can be read at \cite{Schell}, \cite{burgatti} and \cite
{Tric}. Almost half a century ago,  \cite{1}, collected many problems on the base of the constraint type, and all solved by the elliptic integrals
of first and second kind, while the third kind and Theta functions appear 
marginally for $3-D$
deformations.
In order to set the record straight, the concept of \textit{rod} will be recalled from \cite{Goss}:
\begin{quote}
A ‘‘rod’’ is something of a hybrid, i. e. a mathematical curve made up of ‘‘material points’’. It has no cross section, yet it has stiffness. It can weigh nothing or it can be heavy. We can twist it, bend it, stretch it and shake it, but we cannot break it, it is completely elastic. Of course, it is a mathematical object, and does not exist in the physical world; yet it finds wide application in structural (and biological) mechanics: columns, struts, cables, thread, and DNA have all been modeled with rod theory.

\end{quote}
In common practice, a simplified expression to the elastica's curvature is currently used, while the real exact expression leads to nonlinearities.

The link between elliptic functions and elastica was always understood to be close, to
induce the elliptic functions to be plotted through suitable curved rubber rods \cite{green}, while the paper \cite{Basoco} is
founded on the use of the Weierstrass $\wp $ function. Even in contemporary literaure the touchs of the special functions with elastica problems are very close, \cite{1}. In fact, the situations one can meet require either a numerical \cite{wag}, or an elliptic approach \cite{turzi, el1, el2, el3} or a hypergeometric \cite{Hindawi}. In this paper we will follow the last one.

\section{Aim of the paper}

There is a large amount of literature providing the true elastica by solving the relevant nonlinear deflection equation \eqref {fond}, where the rod's curvature is forced by a bending moment. As a matter of fact in  \cite{Hindawi}, we solved the nonlinear deflection problem, as a pure strain analysis, for slender rods by means of Lauricella hypergeometric functions. 

Nevertheless, as far as we are concerned,  the nonlinear, deflection analysis has never been used for a further static analysis until now.

We refer to {\it statically indeterminate} structures, of which the  degree of redundancy is the excess of constraints beyond the strictly necessary. This terminology is currently used as one can see in the literature, see for instance \cite{R1}, or the Wikipedia entry \lq\lq Statically indeterminate''.

Such systems require additional equations and cannot be analyzed through the equilibrium equations alone. On the purpose, a frequent approach, the so--called {\it forces method}, obtains them by transforming the given redundant structure into a so--called {\it simple} one, i.e. a statically determinate one, by  replacing the redundant constraints with the unknown reactions. This will not introduce alterations in stress distribution and in deflections. By decoupling such a multiple load condition over the simple structure in each load stage, each corresponding displacement is computed as a function of the system parameters (say $E,\, J,\, L,\,\ldots$) and the redundant unknowns. Finally, one shall require the consistency, namely that the displacement composition meets the real situation. And so we can  have consistency to a zero displacement (perfect constraint), or to a certain imposed value (anelastic yield) or  to a damper effect (elastic yield): in such a way the necessary elasticity equations are finally obtained. All the computing work one does about such deflections of straight rods of span $L$ is grounded upon the small strains assumption: if $y(x)$ is the deformed centerline's equation, for each $x$ it shall be:
\begin{equation}
y_{x}^2\ll1\quad \text{with $0\leq x\leq L$.}
\label{restriction}
\end{equation}
In such a way, all further steps -and the statically indeterminate reactions- lose their validity when the load intensity, beam slenderness, or high flexibility, causes the relevant strain-even kept in the linear elasticity (effect proportional to its cause) to be such that \eqref {restriction} becomes.
In such cases, the exact centerline curvature is necessary and the nonlinear ODE \eqref {bern} of deflections is often employed in advanced applications regarding the strength of materials in the contexts of aerospace or satellites, see \cite{bele}.

We are just going, on the ground of these new considerations, to compute the statical indeterminate unknowns through a nonlinear deflection analysis.

\section{Statement of the problem}
Let us tackle a $L$-long, thin rod whose end is clamped and A
free. We put at A the origin of a $\left( x,y\right) $ cartesian reference
frame with $x$ along the undeformed rod, from the wall to A; and $y$ downwards,
normal at A to $x$, see \figurename~\ref{f01}. Our basic assumptions are:

\begin{itemize}
\item[A$_1$]  the rod is thin, initially straight, homogeneous, with a uniform
cross section and uniform flexural stiffness $EJ$, where\ $E$ is the Young
modulus, and $J,$ the cross section moment of inertia about a ``neutral''
axis normal to the plane of bending and passing through the central line;

\item[A$_2$]  the slender rod is always charged by coplanar dead (not follower) loads;

\item[A$_3$]  Linear constitutive load holds: so that the induced curvature is
proportional via $1/(EJ)$ to the bending moment.  

\item[A$_{4}$]  The shear transverse deformation is ignored.

\item[A$_{5}$]  A stationary strain field, by isostatic equilibrium of
active loads and reactive forces takes place, and, due to rest of static
equilibrium, no rod element undergoes acceleration.
\end{itemize}

The key Bernoulli-Euler linear constitutive law connects flexural stiffness $%
EJ$ and bending moment $\mathscr{M}$ to the consequent curvature $\varkappa $ of the
deflected centerline: $\pm EJ\varkappa =\mathscr{M}(x).$ With our (traditional)
reference frame $(x,y)$, the curvature always bears the opposite sign to the
bending moment: 
\begin{equation}
EJ\frac{y_{xx}}{\left( 1+y_{x}^{2}\right) ^{3/2}}=-\mathscr{M}(x)  \label{bern}
\end{equation}
where $y=y(x)$ is the unknown elastica's equation and 
\[
y_{x}=\dfrac{\mathrm{d}y}{\mathrm{d}x},\quad y_{xx}=\dfrac{\mathrm{d^{2}}y}{%
\mathrm{d}x^{2}}. 
\]
To (\ref{bern}) a useful dress can be given, so that it shall be minded that the Cauchy problem: 
\begin{equation}
\begin{cases} \dfrac{{\rm d}}{{\rm
d}x}\left(\dfrac{y_{x}}{\sqrt{1+y_x^{2}}}\right)=-
\dfrac{1}{{EJ}}\,\mathscr{M}(x),\quad 0\leq x\leq L \\ y(L)=0\\ y_x(L)=0 \end{cases}
\label{fond}
\end{equation}
can be solved in such a way. Putting: 
\[
H(x):=\int_{x}^{L}\mathscr{M}(\xi )\mathrm{d}\xi 
\]
the solution of \eqref{fond} is: 
\begin{equation}
y(x)=-\int_{x}^{L}\frac{H(\xi )}{\sqrt{\left( {EJ}\right) ^{2}-H^{2}(\xi )}}%
\,\mathrm{d}\xi .  \label{tesi}
\end{equation}
In order to have a solution defined for $x\in\left[0,L\right]$ we have to require that
\begin{equation}\label{globbal}
-1<\frac{1}{EJ}H(x)<1
\end{equation}
which is a general prescription on the load effects and will be defined in single cases. 
\subsection{Nonlinear deflection, a first example}
For a $L,\, E,\, J$  cantilever tip-sheared by $P$, the approximate linear theory provides a tip's vertical displacement $\delta_{0}=(PL^3)/(3EJ)$. In such a case, in \eqref{fond} we put $\mathscr{M}(x)=-Px$ and defining the nondimensional quantities:
$$
\mu=\frac{PL^3}{3EJ}; \quad \eta=\frac{y}{L};\quad \xi=\frac{x}{L}
$$
it becomes $\delta_{0}=2\mu/3$. 

Formula (2.32) of our previous paper \cite{Hindawi} for the exact free-end displacement provides:
\[
\eta(\xi)=\frac{\mu}{3 \left(1-\mu ^2\right)^{3/2}}\, \xi ^3-\frac{\mu 
   }{\sqrt{1-\mu ^2}}\, \xi+ \sqrt{\frac{\mu}{2}}\,A(\mu)
\]
where $A(\mu)$ is a complicated function involving elliptic integrals of first and third kind according to formula (2.29) of the same \cite{Hindawi}. They can be expanded with respect to $\mu$, with initial point $\mu_0=0$; then, recalling the Maclaurin first--order expansion with respect to $\mu$
\begin{equation}\label{f2}
\frac{\mu}{3 \left(1-\mu ^2\right)^{3/2}}\, \xi ^3-\frac{\mu 
   }{\sqrt{1-\mu ^2}}\, \xi=\mu  \left(\frac{\xi ^3}{3}-\xi \right)+O\left(\mu ^2\right)
\end{equation}
one obtains eventually
\begin{equation}\label{bell}
\eta(\xi)=\frac{\mu}{3}\left(2-3\xi+\xi^3\right)+O\left(\mu ^2\right)
\end{equation}
i.e. the usual approximation which for $\xi=0$ (free tip) provides the same deflection
$\delta_{0}$ as above.

\begin{figure}
\begin{center}
\includegraphics{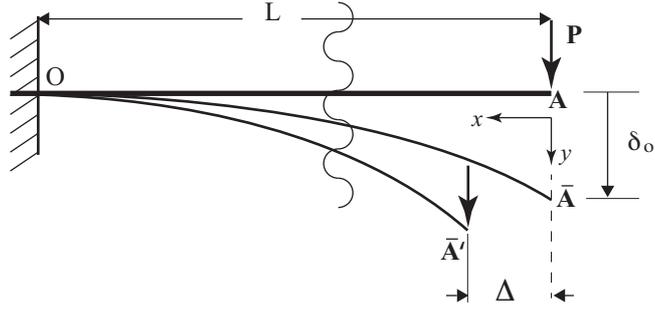}
\caption{Cantilever's kinds of deformed centerline and transverse displacement $\Delta$}\label{f01}
\end{center}
\end{figure}

\figurename~\ref{f01} shows qualitatively the behavioral differences, for a $L$-rod $OA$ under a dead load $P$ at the tip, between its typical approximate elastica  $O\overline{A}$ and the true elastica $O\overline{A'}$. As a matter of fact, the exact treatment highlights a horizontal displacement of the free tip which is given by
\[
\Delta =L-\inf_{0\leq y\leq x(L)}\left[x(y)\right]
\]
and which has been computed (sect. 2.1 and 3.) in \cite{Hindawi}.

What we are now going to do is to employ a displacement computation like this for gaining exact elasticity equations to solve statically indeterminate systems.

\section{The statically indeterminate heavy cantilever supported by a roller}

\begin{figure}
\begin{center}
\scalebox{1.0}{\includegraphics{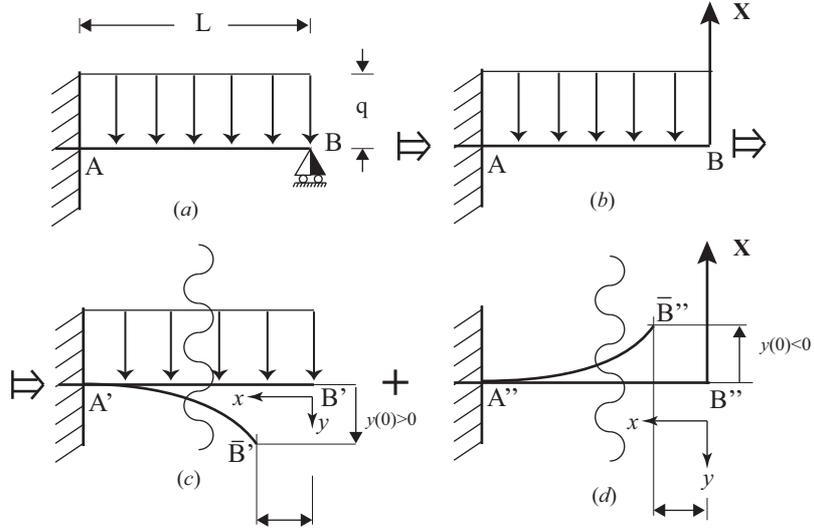}}
\caption{The heavy cantilever supported by a roller (a), the relevant {\it simple} structure (b) and the  loads decoupling (c) and (d).}\label{f02}
\end{center}
\end{figure}

Our first statically indeterminate system is a rod (\figurename~\ref{f02} (a)) under a uniform load $q$, whose one side A is built-in, while the free tip B is supported by a roller: for such a 1-redundant rod we assume as unknown the reaction $\mathscr{X}$ just acted by the roller. Then the statically \lq \lq simple'' rod is of (\figurename~\ref{f02} (b)), whose analysis can be decoupled according to \figurename~\ref{f02} (c) and (d) and tackled separately.

\subsection{First sub-system: the heavy cantilever}
Let us consider first the subsystem of \figurename~\ref{f02} (c). For the bending moment we have $\mathscr{M}(x)=-(q/2)x^2.$ Specializing \eqref {globbal} for this load, it will imply $q<6EJ/L^3.$ Thus \eqref{tesi} reads\footnote{The approach of computing the exact tip position of the rod would be  extremely long and capable of hiding  the sense of this research deviating by its central purpose. In such a way the (very small) horizontal deflection of the tip has been put to zero. Anyway, an example where the rod is assumed {\it inextensible} and the axial deformation is kept in account, can be seen in our previous paper [14] sect. 3, formula (3.12).
} as
\begin{equation}\label{solivp1}
y(x)=\int_x^L\frac{q(L^3-y^3)}{6\sqrt{1-\left(\frac{EJ}{6}\right)^2(L^3-y^3)^2}}\,{\rm d}y
=\frac{2EJ}{qL^2}\int_0^{\frac{q}{6EJ}(L^3-x^3)}\frac{u}{(1-u^2)^{1/2}\left(1-\frac{6EJ}{qL^3}u\right)^{2/3}}\,{\rm d}u
\end{equation}
Now using the Lauricella functions' integral representation theorem, we infer the integration formula
\begin{equation}\label{irtfd3a}
\int_0^a\frac{u}{(1-u^2)^{1/2}(1-bu)^{2/3}}{\rm d}u=\frac{a^2}{2}\mathrm{F}_{D}^{(3)}\left( \left. 
\begin{array}{c}
2;1/2,1/2,2/3 \\[2mm]
2
\end{array}
\right| a,-a,ab\right) 
\end{equation}
Clearly \eqref{irtfd3a} is the key to evaluate the integral in \eqref{solivp1}, so that
\begin{equation}\label{solivp1h}\tag{\ref{solivp1}b}
y(x)=\frac{q \left(L^3-x^3\right)^2}{36 E J L^2}\mathrm{F}_{D}^{(3)}\left( \left. 
\begin{array}{c}
2;1/2,1/2,2/3 \\[2mm]
3
\end{array}
\right| \frac{q}{6EJ}(L^3-x^3),-\frac{q}{6EJ}(L^3-x^3),1-\frac{x^3}{L^3}\right) 
\end{equation}
which provides the deflection law for each $x$. We are interested to take $x=0$ in \eqref{solivp1h} getting
\begin{equation}\label{solivp1k}\tag{\ref{solivp1}c}
y(0)=\frac{L^4 q}{36EJ}\mathrm{F}_{D}^{(3)}\left( \left. 
\begin{array}{c}
2;1/2,1/2,2/3 \\[2mm]
3
\end{array}
\right|\frac{L^3 q}{6EJ},-\frac{L^3 q}{6EJ},1\right) 
\end{equation}
Notice that for $x=0$ the third argument of Lauricella's becomes 1, so that such a function collapses into an Appell one and from there it can be reduced to a $_3\mathrm{F}_{2}$ hypergeometric function, as shown in section \ref{identities}. In such a way we will be allowed to describe the strain by means of a function of only one variable. We can simplify \eqref{solivp1k} using the reduction formulae \eqref{toeta0b} and \eqref{wolfram} in sequence,
getting
\begin{equation}\label{solivp1y}\tag{\ref{solivp1}d}
y(0)=\frac{L^4 q}{8EJ}\,_3\mathrm{F}_{2}\left( \left. 
\begin{array}{c}
\frac{1}{2},1,\frac32\\[2mm]
\frac{7}{6},\frac53
\end{array}
\right|\frac{L^6 q^2}{36(EJ)^2}\right)
\end{equation}
which provides as a function of $q,\, L,\, E,\, J$  the free tip true deflection of rod of \figurename~\ref{f02} (c).
\subsection{Second sub-system: the tip-sheared cantilever}
Allow us consider the subsystem of \figurename~\ref{f02} (d). The bending moment being $\mathscr{M}(x)=\mathscr{X}x,$ so that \eqref{tesi} reads as
\begin{equation}\label{solivp2}
y(x)=\int_x^L\frac{\mathscr{X}(y^2-L^2)}{\sqrt{(2EJ)^2-\mathscr{X}^2(y^2-L^2)^2}}\,{\rm d}y=\frac{\mathscr{X}}{4EJL}\int_{x^2-L^2}^{0}\frac{u}{\left(1-\left(\frac{\mathscr{X}}{2EL}u\right)^2\right)^{1/2}\left(1+\frac{1}{L^2}u\right)^{1/2}}\,{\rm d}u
\end{equation}
With the integral representation theorem we also infer the integration formula
\begin{equation}\label{irtfd3b}
\int_0^a\frac{u}{(1-b^2u^2)^{1/2}(1+cu)^{1/2}}{\rm d}u=\frac{a^2}{2}\mathrm{F}_{D}^{(3)}\left( \left. 
\begin{array}{c}
2;1/2,1/2,1/2 \\[2mm]
3
\end{array}
\right| ab,-ab,-ac\right) 
\end{equation}
Clearly, \eqref{irtfd3b} is the key to evaluate the integral in \eqref{solivp2}, giving
\begin{equation}\label{solivp2h}\tag{\ref{solivp2}b}
y(x)=-\frac{\left(x^2-L^2\right)^2\mathscr{X}}{8EJL}\mathrm{F}_{D}^{(3)}\left( \left. 
\begin{array}{c}
2;1/2,1/2,1/2 \\[2mm]
3
\end{array}
\right| \frac{(x^2-L^2)\mathscr{X}}{2EJ},-\frac{(x^2-L^2)\mathscr{X}}{2EJ},1-\frac{x^2}{L^2}\right) 
\end{equation}
We are interested to take $x=0$ in \eqref{solivp2h}, getting
\begin{equation}\label{solivp2k}\tag{\ref{solivp2}c}
y(0)=-\frac{L^4 \mathscr{X}}{8EJ}\mathrm{F}_{D}^{(3)}\left( \left. 
\begin{array}{c}
2;1/2,1/2,1/2 \\[2mm]
3
\end{array}
\right|-\frac{L^2\mathscr{X}}{2EJ},\frac{L^2\mathscr{X}}{2EJ},1\right) 
\end{equation}
We can simplify \eqref{solivp2k} using the reduction formulae \eqref{toeta0b} and \eqref{wolfram} in sequence,
getting
\begin{equation}\label{solivp2y}\tag{\ref{solivp2}d}
y(0)=-\frac{L^3 \mathscr{X}}{3EJ}\,_3\mathrm{F}_{2}\left( \left. 
\begin{array}{c}
\frac{1}{2},1,\frac32\\[2mm]
\frac{5}{4},\frac74
\end{array}
\right|\frac{L^4 \mathscr{X}^2}{4(EJ)^2}\right)
\end{equation}

\subsection{Consistency}
In order to evaluate the redundant unknown $\mathscr{X}$, notice that the addition of two displacements given by \eqref{solivp1y} and \eqref{solivp2y} shall be zero, arriving at the transcendental equation:
\begin{equation}\label{trans}
3Lq\, \,_3\mathrm{F}_{2}\left( \left. 
\begin{array}{c}
\frac{1}{2},1,\frac32\\[2mm]
\frac{7}{6},\frac53
\end{array}
\right|\frac{L^6 q^2}{36(EJ)^2}\right)-8\mathscr{X}\,_3\mathrm{F}_{2}\left( \left. 
\begin{array}{c}
\frac{1}{2},1,\frac32\\[2mm]
\frac{5}{4},\frac74
\end{array}
\right|\frac{L^4 \mathscr{X}^2}{4(EJ)^2}\right)=0
\end{equation}
in our unknown reaction  $\mathscr{X}$. Of course, \eqref{trans} cannot be solved analytically, but numerically.
Nevertheless we will provide an analytical approximate solution. In fact, if we approximate both hypergeometrical series only by their first term which is equal to 1, the solution to \eqref{trans} is approximately given by:
\[
\mathscr{X}=\frac38\,Lq
\]
according to the linearized theory.

For a higher approximation, we can proceed through the Lagrange series reversion.\footnote{Lagrange, J.L. (1770): "Nouvelle méthode pour résoudre les équations littérales par le moyen des séries," Mémoires de l'Académie Royale des Sciences et Belles-Lettres de Berlin, vol. 24, pages 251–326,  a paper of 1766 but issued four years later. The Lagrange {\it inversion} theorem used here gives the Taylor series expansion of the inverse function of an analytic function; it shall not be confused with Lagrange {\it reversion} theorem which provides formal power series expansions of certain functions defined implicitly.}. 

For the purpose, in \eqref{trans} we put:
\[
\mathscr{X}=\frac{2EJ}{L^2}\mathscr{Y}
\]
so that \eqref{trans} can be written as
\begin{equation}\tag{\ref{trans}b}
f(\mathscr{Y}):=\mathscr{Y}\,_3\mathrm{F}_{2}\left( \left. 
\begin{array}{c}
\frac{1}{2},1,\frac32\\[2mm]
\frac{5}{4},\frac74
\end{array}
\right|\mathscr{Y}^2\right)=\frac{3}{16}\frac{L^3q}{EJ}\,_3\mathrm{F}_{2}\left( \left. 
\begin{array}{c}
\frac{1}{2},1,\frac32\\[2mm]
\frac{7}{6},\frac53
\end{array}
\right|\frac{L^6 q^2}{36(EJ)^2}\right)=\mathscr{Z}
\end{equation}
Function $f:]-1,1[\to\mathbb{R}$ appearing at left hand side of (\ref{trans}b) is odd, bijective and then invertible. Notice that, due to the hypergeometric series divergent behaviour at the edge of the convergence disk, we have: $f(\mathscr{Y})\to+\infty$ as $\mathscr{Y}\to1^{-}.$

\begin{figure}
\begin{center}
\scalebox{0.8}{\includegraphics{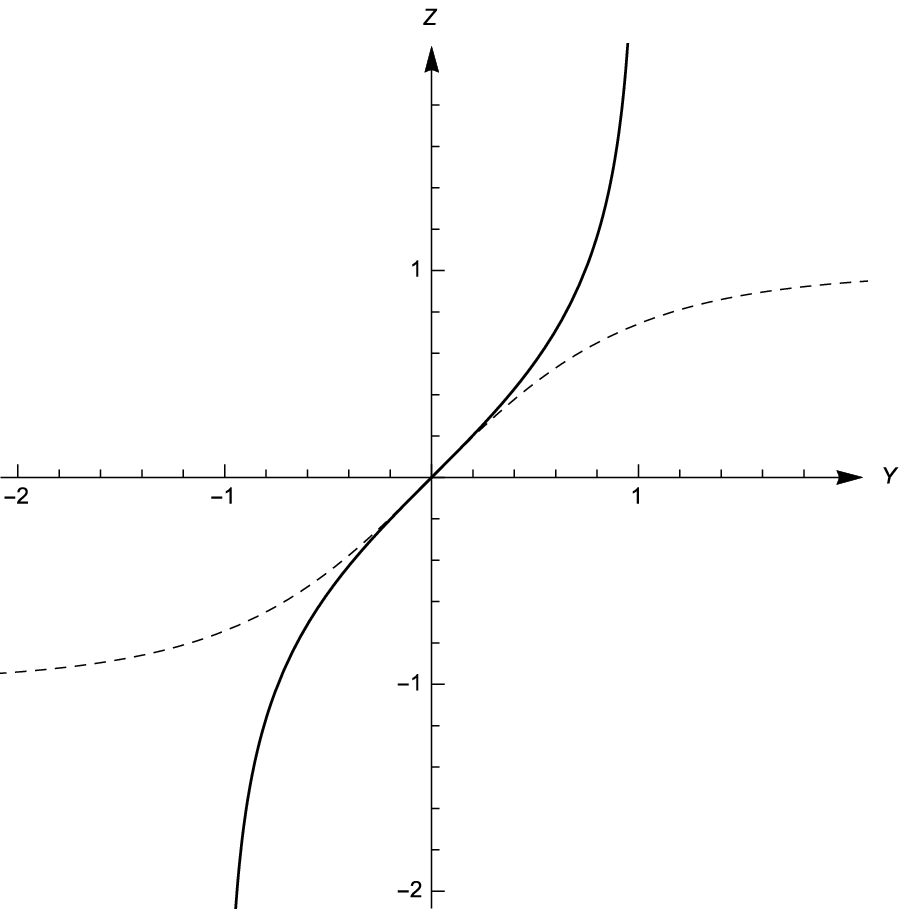}}
\caption{$f(\mathscr{Y})$ and $f^{-1}(\mathscr{Z}) $}\label{f03}
\end{center}
\end{figure}

\begin{figure}
\begin{center}
\boxed{\scalebox{0.75}{\includegraphics{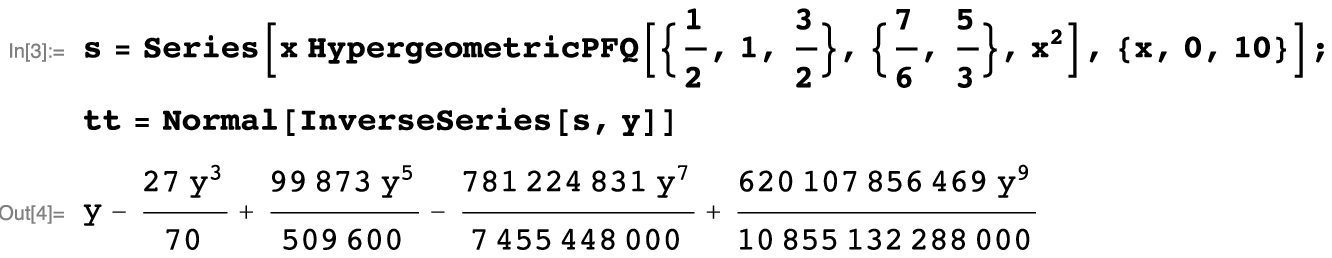}}}
\caption{Series for $\mathscr{Z}=f^{-1}(\mathscr{Y})$}\label{M01}
\end{center}
\end{figure}

In \figurename~\ref{f03} we depict continuously the plot of $f(\mathscr{Y})$ and as a dotted line the inverse $f^{-1}(\mathscr{Z}).$
By Mathematica, based on the Lagrange {\it inversion} theorem, we can obtain the series expansion with respect to $\mathscr{Z}$ of the root of:  $f(\mathscr{Y})=\mathscr{Z}$. We give hereinafter the relevant computer instructions in \figurename~\ref{M01}

After this, we will expand the right hand side of (\ref{trans}b),  evaluate the approximate solution and expand it in power series:
\begin{figure}
\begin{center}
\boxed{\scalebox{0.75}{\includegraphics{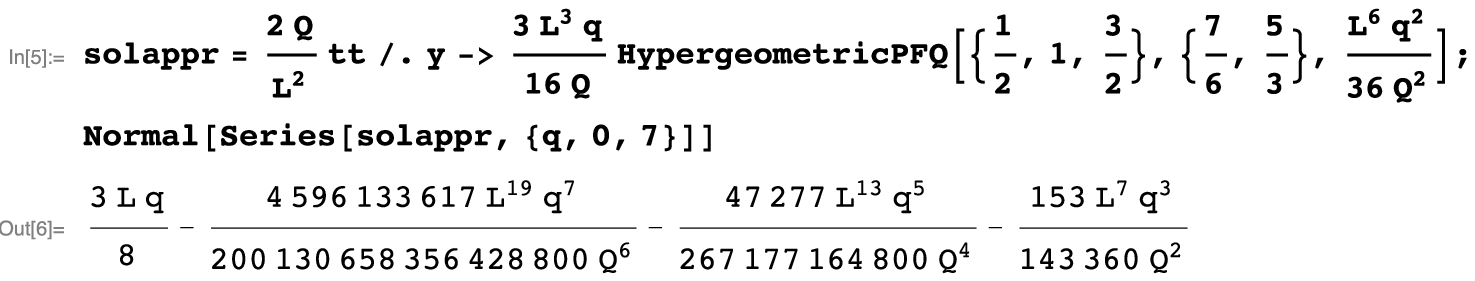}}}
\caption{Approximate solution of (\ref{trans}b)}\label{M02}
\end{center}
\end{figure}

\begin{figure}
\begin{center}
\scalebox{0.8}{\includegraphics{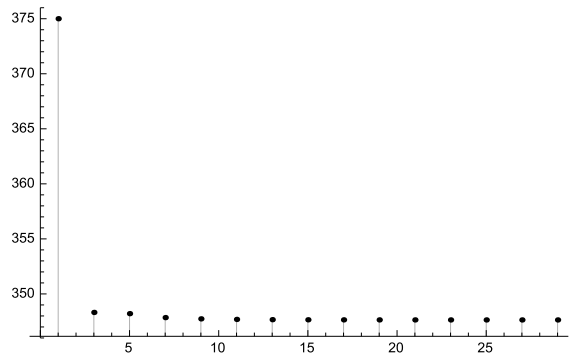}}
\caption{The statically indeterminate heavy cantilever supported by a roller with $L=1,\,q=1000,\,EJ=200$, $\mathscr{X}$ values versus number of terms of its expansion}\label{f04}
\end{center}
\end{figure}

Now, going back to our unknown $\mathscr{X}$, we can get an approximate solution of (\ref{trans}b):
  \begin{equation}
\mathscr{X}=\frac{3 L q}{8}-\frac{153 L^7 q^3}{143360 (EJ)^2}-\frac{47277 L^{13} q^5}{267177164800(EJ)^4}-\frac{4596133617 L^{19} q^7}{200130658356428800(EJ)^6}+\cdots
\label{transaa} 
\end{equation}
The plot in \figurename~\ref{f04} shows the behavior of $\mathscr{X} $ versus the number of terms of  \eqref{transaa} wherefrom the linearized solution is coming for $n=0$ and overestimates the solution; the addition of further terms provides us with a better approximation which is stable just after $n=7$.  
\section{The statically indeterminate two-side built-in heavy rod}
\begin{figure}
\begin{center}
\scalebox{1.0}{\includegraphics{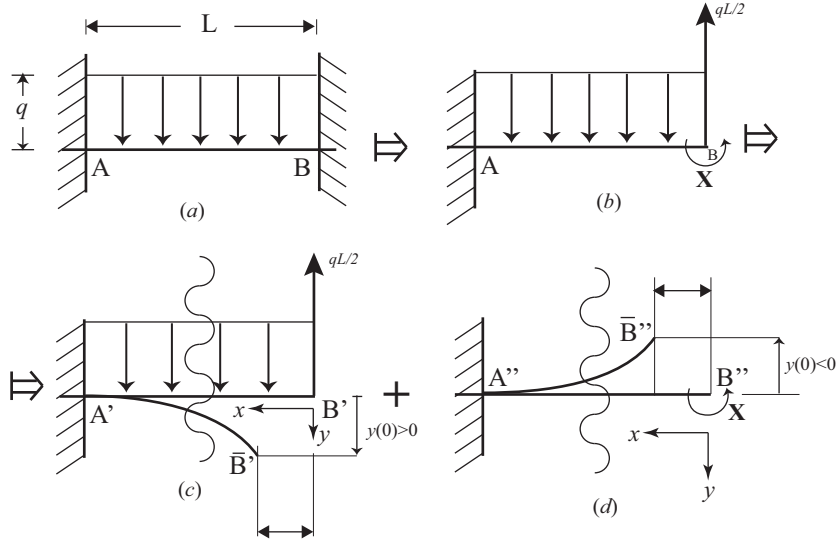}}
\caption{The doubly built-in cantilever (a), the  relevant {\it simple} structure (b) and the loads decoupling (c) and (d). }\label{f05}
\end{center}
\end{figure}

Our second problem, \figurename~\ref{f05}, will concern the constraints' redundancy of a doubly built-in slender rod AB subjected to a uniform continuous load $q$. Such a system, for a general load, is three times redundant; for a purely vertical but non--symmetrical one, two times. In our case the actual degree will be only one. We assume as unknown $\mathscr{X}$ the bending moment acted by one of the built-in supports (see \figurename~\ref{f05} (b)). The same constraint will provide a vertical reaction equilibrating just half a load and then of value $qL/2$. In such a way, the bending moment at the $x$-section will be given by:
\[
\mathscr{M}(x)=\frac{qLx}{2}-\frac{x^2}{2},\quad 0\leq x\leq L .
\] 
As before, we are going to analyze apart the subsystems (c) and (d).

\subsection{First sub-system: a heavy tip-sheared cantilever}
Let
\[
\mathscr{M}(x)=\frac{L q x}{2}-\frac{q x^2}{2}
\]
for the rod of \figurename~\ref{f05} (c). Observe that \eqref{globbal} reads as
\[
-1<-\frac{1}{12EJ} q (L-x)^2 (L+2 x)<1
\]
which will imply the condition $q<12EJ/L^3$ to be met.

Accordingly, \eqref{tesi} leads\footnote{Being the double clamping made for guarantee two fixed positions, the maximum range of integration is from $0$ to $L.$ Then we fix $x\in[0,L]$, which implies the deformed rod is {\it extensible} under its load.
} to the hyperelliptic integral:
\begin{equation}\label{solivpbis}
y(x)=-\int_x^L\frac{q (L-y)^2 (L+2 y)}{\sqrt{144 (EJ)^2-q^2 (L-y)^4 (L+2 y)^2}}{\rm d}y
\end{equation}
and in particular, for the free tip, obviously:
\begin{equation}\label{solivp1z}\tag{\ref{solivpbis}b}
y(0)=-\int_0^L\frac{q (L-y)^2 (L+2 y)}{\sqrt{144 (EJ)^2-q^2 (L-y)^4 (L+2 y)^2}}{\rm d}y:=-\mathcal{I}(L,q,E,J)
\end{equation}

\subsection{Second sub-system: a cantilever under a bending moment at its free end}

Let $\mathscr{M}(x)=\mathscr{X}$ for the rod of \figurename~\ref{f05} (d). In such a way \eqref{tesi} is solved by elementary methods:
\begin{equation}\label{solivpbish}
y(x)=\frac{\sqrt{(EJ)^2-\mathscr{X}^2 (x-L)^2}-EJ}{\mathscr{X}}
\end{equation}
which gives in particular
\begin{equation}\label{solivp2z}\tag{\ref{solivpbish}b}
y(0)=\frac{\sqrt{(EJ)^2-\mathscr{X}^2L^2}-EJ}{\mathscr{X}}
\end{equation}
and then we succeed in evaluating the displacement provided by the redundant unknown $\mathscr{X}$ if it would act alone.
\subsection{Consistency}
\begin{figure}
\begin{center}
\scalebox{0.8}{\includegraphics{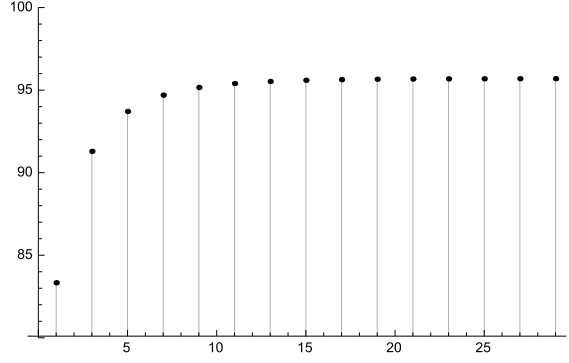}}
\caption{The doubly built-in cantilever with $L=1,\,q=1000,\,EJ=200$: $\mathscr{X}$ values versus number of terms of its expansion}\label{f06}
\end{center}
\end{figure}

To compute the redundant unknown $\mathscr{X}$, our constraint being perfect, we shall require again the addition of two displacements given by \eqref{solivp1z} and \eqref{solivp2z} shall be zero. The computation of the integral \eqref{solivp1z}, because of the cubic, would lead to Lauricella higher order functions. By expanding in powers the integrand in \eqref{solivp1z}, we get the approximation:
\begin{equation}\label{appproxx}
\mathcal{I}(L,q,E,J)\simeq-\frac{L^4q}{24EJ}\,_2\mathrm{F}_{1}\left( \left. 
\begin{array}{c}
\frac{1}{2},\frac23\\[2mm]
\frac{5}{3}
\end{array}
\right|\frac{L^6 q^2}{36 (EJ)^2}\right)
\end{equation}
Then, solving to $\mathscr{X}$ the displacement equation:
\[
\frac{\sqrt{(EJ)^2-\mathscr{X}^2L^2}-EJ}{\mathscr{X}}=\mathcal{I}(L,q,E,J)\implies \mathscr{X}=\frac{2 EJ \mathcal{I}(L,q,E,J)}{L^2+\mathcal{I}^2(L,q,E,J)} 
\]
Putting to  $\mathcal{I}(L, q, E, J )$ its approximation given by \eqref{appproxx} and after a further development, we get:
\begin{equation}\label{transbb}
\mathscr{X}=\frac{L^2 q}{12}+\frac{11 L^8 q^3}{34560 (EJ)^2}+\frac{77 L^{14} q^5}{19906560 (EJ)^4}+\cdots
\end{equation}
 which holds even powers of $EJ.$
\figurename~\ref{f06} shows the behavior of such a solution with the number of terms of the series \eqref{transbb}: the linearized solution is that for $n=0$ and underestimates the  $\mathscr{X}$; the addition of further terms provides us with a better approximation which is stable after $n=11$. 

\section{Conclusions}
In order to get a better clarity, we will refer to a sample problem.
Let us take into account a slender rod  of length $L=1 m$ and a square uniform section of side $s=10^{-2}m$, so that $J=s^{4}/12 =10^{-9} m^{4}$. The material is steel with $E = 2.1\cdot10^{11}N/m^{2}$, so that each section has a flexural stiffness $EJ\simeq  200  N m^{2}$. 

Let us consider first the statically indeterminate unknown relevant to \figurename~\ref{f02} : according to the approximate linear theory, it is given by $3qL/8$, namely $375 N$ for our sample problem.
We take a unit continuous load of $q=10^{3} N/m$ which is compliant to the prescription of being less than $6EJ/L^{3}$.

Our nonlinear approximate solution is provided by a relationship of kind \eqref{transaa}, namely:
$$
\mathscr{X}=\sum_{n=0}^\infty a_{n}\frac{L^{6n+1}q^{2n+1}}{(EJ)^{2n}}
$$
The reaction $\mathscr{X}$ of the statically indeterminate roller is then computed via a power series involving length, flexural rigidity, unit load and the relevant coefficients, some of which are shown by \eqref{transaa}. Therefore \figurename~\ref{f04} shows that the linearized approach overestimates about 7.8\% the \lq\lq reference value'' provided by the (less approximate) nonlinear treatment which converges from $n=7$ on.

Passing to the second statically indeterminate system (see \figurename~\ref{f05}) we took the same unit continuous load $q=10^{3} N/m$ which is compliant to the new prescription of being less than $12 EJ/L^{3}$. According to the approximate linear theory, the statically indeterminate bending moment (hereinafter called again $\mathscr{X}$ ) exerted by one of the building-in constraints, is traditionally accepted as $qL^{2}/12$, namely $83.33Nm$ for our sample problem. But a proper number of terms  of  \eqref{transbb} 
$$
\mathscr{X}=\sum_{n=0}^\infty b_{n}\frac{L^{6n+2}q^{2n+1}}{(EJ)^{2n}}, 
$$
provides a value of $95.16 Nm$ which is 14.19\% greater.

The general conclusions suggested by all of the above are the following:

\begin{itemize}

\item[C$_1$]  
In both problems analyzed in this paper, the traditional linearized approach leads to significant errors on the statically unknowns. Such deviations increase for rods that are progressively more and more slender and flexible, see \figurename~\ref{f04}, \figurename~\ref{f06}.

\item[C$_2$]
The above mentioned effects obviously lead to a remarkable impact on the axial bending stresses, say $\sigma$ s. E.g. for the rod in \figurename~\ref{f02}, the highest bending moment is given by $\mathscr{X}^{2}/2q $ and in that section the relevant edge $\sigma$ is  16\% lower than that from the linearized approach.
\item[C$_{3}$]  
For instance, referring to the rod in \figurename~\ref{f05},  in the section of the highest moment, the relevant edge $\sigma$ is 6\% greater than that from the linearized approach. 
\item[C$_4$]  It is commonly believed and written that the redundant unknowns do not depend on the flexural stiffness, but this is, generally speaking, not true.

\end{itemize}

In fact, our final formulae \eqref{transaa} and \eqref{transbb} show undoubtlessly that this is due to having taken $n=0$, namely only the first term of the expansion whose subsequent terms hold the reciprocal of powers of $EJ$.
Therefore the Belluzzi sentence, see \cite{Bellu} p. 363:
\begin{quote}

Si osservi che quando i vincoli sono rigidi la reazione $A$ non dipende da $EJ$, mentre ne dipende quando sono cedevoli.\footnote{Notice that when the constraints are rigid, the (redundant) reaction does not depend upon $EJ$, whilst it does when they undergo some yield.}

\end{quote}
can be accepted {\it solely within a first order approximation}.

The linearized approach does its approximation at the very beginning of the treatment, so that leads to a less refined and elementary conclusion. On the contrary, we keep the nonlinear nature of the problem, generating a true solution which -even if we decide/need to approximate it at final stages- reveals new mathematical links and conclusions.


\section*{Appendix: Hypergeometric identities}\label{identities}

We recall hereinafter something about the non too known Lauricella
hypergeometric functions.

The first hypergeometric series appeared in the Wallis's \textit{Arithmetica
infinitorum} (1656): 
\[
_{2}F_{1}(a,b;c;x)=1+\frac{a\cdot b}{1\cdot c}x+\frac{a\cdot (a+1)\cdot b%
\cdot (b+1)}{1\cdot 2\cdot c\cdot (c+1)}x^{2}+\cdots ,
\]
for $|x|<1$ and real parameters $a,\,b,\,c.$ The product of $n$ factors: 
\[
(\lambda )_{n}=\lambda \left( \lambda +1\right) \cdots \left( \lambda
+n-1\right) ,
\]
\noindent called \textit{Pochhammer symbol} (or \textit{truncated factorial}
) allows to write $_{2}F_{1}$ as: 
\[
_{2}F_{1}(a,b;c;x)=\sum_{n=0}^{\infty }\frac{\left( a\right) _{n}\left(
b\right) _{n}}{\left( c\right) _{n}}\frac{x^{n}}{n!}.
\]
A meaningful contribution on various $_{2}F_{1}$ topics
is ascribed to Euler\footnote{%
We quote three works: a) \textit{De progressionibus transcendentibus}, Op.
omnia, S.1, vol. 28; b) \textit{De curva hypergeometrica} Op. omnia, S.1,
vol. 16; c) \textit{Institutiones Calculi integralis}, 1769, vol. II};
but he does not seem \cite{dutka} to have known the integral representation: 
\[
_{2}F_{1}(a,b;c;x)=\frac{\Gamma (c)}{\Gamma (a)\Gamma (c-a)}\,\int_{0}^{1}%
\frac{u^{a-1}(1-u)^{c-a-1}}{(1-xu)^{b}}\,\mathrm{d}u,
\]
 really due to A. M. Legendre\footnote{%
A. M. Legendre, \textit{Exercices de calcul int\'{e}gral}, II, quatri\'{e}me
part, sect. 2, Paris 1811}. The above integral relationship is true if $c>a>0$ and for 
$\left| x\right| <1,$ even if this limitation can be discarded thanks to the analytic
continuation.\ 

Furthermore, many functions have been introduced in 19$^{\mathrm{th}}$ century for
generalizing the hypergeometric functions to multiple variables. We recall the Appell ${\rm F}_1$ two--variable hypergeometric series, defined as: 
\[
\mathrm{F}_{1}\left( \left. 
\begin{array}{c}
a;b_{1},b_{2} \\[2mm]
c
\end{array}
\right| x_{1},x_{2}\right) =\sum_{m_{1}=0}^{\infty }\sum_{m_{2}=0}^{\infty }%
\frac{(a)_{m_{1}+m_{2}}(b_{1})_{m_{1}}(b_{2})_{m_{2}}}{(c)_{m_{1}+m_{2}}}%
\frac{x_{1}^{m_{1}}}{m_{1}!}\frac{x_{2}^{m_{2}}}{m_{2}!},\quad |x_{1}|<1,\,|x_{2}|<1
\]
The analytic continuation of  Appell's function on
$\mathbb{C}\setminus [1,\infty)\times\mathbb{C}\setminus [1,\infty)$ comes from its integral representation theorem: if $\rm{Re}a>0,\,\rm{Re}(c-a)>0$: 
\begin{equation}
\mathrm{F}_{1}\left( \left. 
\begin{array}{c}
a;b_{1},b_{2} \\[2mm]
c
\end{array}
\right| x_{1},x_{2}\right) =\frac{\Gamma (c)}{\Gamma (a)\Gamma (c-a)}%
\int_{0}^{1}\frac{u^{a-1}\left( 1-u\right) ^{c-a-1}}{\left(
1-x_{1}\,u\right) ^{b_{1}}\left( 1-x_{2}\,u\right) ^{b_{2}}}\,\mathrm{d}u.
\label{F1}
\end{equation}
For us, the functions introduced and investigated by G. Lauricella (1893) and S. Saran (1954), are of prevailing interest; and among them the hypergeometric
function $F_{D}^{(n)}$ of $n\in \mathbb{N}^{+}$ variables (and $n+2$
parameters), see \cite{s} and \cite{l}, defined as: 
\[
F_{D}^{(n)}\left( a,b_{1},\ldots ,b_{n};c;x_{1},\ldots ,x_{n}\right) :=
\sum_{m_{1},\ldots ,m_{n}\in \mathbb{N}}\frac{(a)_{m_{1}+\cdots
+m_{n}}(b_{1})_{m_{1}}\cdots (b_{n})_{m_{n}}}{(c)_{m_{1}+\cdots
+m_{n}}m_{1}!\cdots m_{n}!}\,x_{1}^{m_{1}}\cdots x_{n}^{m_{m}} 
\]
with the hypergeometric series usual convergence requirements $%
|x_{1}|<1,\ldots ,|x_{n}|<1$. 

If $\mathrm{Re}\,c>\mathrm{Re}\,a>0$ , the
relevant Integral Representation Theorem (IRT) provides: 
\[
F_{D}^{(n)}\left( a,b_{1},\ldots ,b_{n};c;x_{1},\ldots ,x_{n}\right) =\frac{%
\Gamma (c)}{\Gamma (a)\,\Gamma (c-a)}\,\int_{0}^{1}\,\frac{%
u^{a-1}(1-u)^{c-a-1}}{(1-x_{1}u)^{b_{1}}\cdots (1-x_{n}u)^{b_{n}}}\,\mathrm{d%
}u 
\]
allowing the analytic continuation to $\mathbb{C}^{n}$ deprived of the
cartesian $n$-dimensional product of the interval $]1,\infty [$ with itself.

Finally, we provide here some reduction formulae employed throughout this paper. 

\begin{equation}\label{r1}
\mathrm{F}_{1}\left( \left. 
\begin{array}{c}
a;b,b,\\[2mm]
a+1
\end{array}
\right|x,-x\right)=\,_2\mathrm{F}_{1}\left( \left. 
\begin{array}{c}
\frac{a}{2};b,\\[2mm]
\frac{a}{2}+1
\end{array}
\right|x^2\right)
\end{equation}
Identity \eqref{r1} follows from the integral representation theorem.
In fact, if we return from the symbolic expression for the Appell function to its meaning in terms of integral representation, we see that
\[
\begin{split}
\mathrm{F}_{1}\left( \left. 
\begin{array}{c}
a;b,b\\[2mm]
a+1
\end{array}
\right|x ,x\right) &=\frac{\Gamma(a+1)}{\Gamma(a)\Gamma(1)}\int_0^1\frac{u^{a-1}}{(1-x^2u^2)^{1/2}}{\rm d}u\\
&=a\int_0^1\frac{v^{\frac{a}{2}-1}}{(1-x^2v)^{1/2}}\,{\rm d}v=\,_2\mathrm{F}_{1}\left( \left. 
\begin{array}{c}
\frac{a}{2};b,\\[2mm]
\frac{a}{2}+1
\end{array}
\right|x^2\right)
\end{split}
\]
Next, we consider some generalizations of the well--known Gauss summation formula: if $c>a+b$ then
\[
\,_{2}\mathrm{F}_{1}\left( \left. 
\begin{array}{c}
a;b \\[2mm]
c
\end{array}
\right| 1\right)=\frac{\Gamma(c)\Gamma(c-a-b)}{\Gamma(c-a)\Gamma(c-b)}
\]
moreover, we also have
\[
\mathrm{F}_{1}\left( \left. 
\begin{array}{c}
a;b_1,b_2 \\[2mm]
c
\end{array}
\right| 1,x\right)=\,_{2}\mathrm{F}_{1}\left( \left. 
\begin{array}{c}
a;b_1 \\[2mm]
c
\end{array}
\right| 1\right)\,\,_{2}\mathrm{F}_{1}\left( \left. 
\begin{array}{c}
a;b_2 \\[2mm]
c-b_1
\end{array}
\right| x\right)
\]
which, when $c>a+b_1$ gives, due to the Gauss summation theorem,
\[
\mathrm{F}_{1}\left( \left. 
\begin{array}{c}
a;b_1,b_2 \\[2mm]
c
\end{array}
\right| 1,x\right)=\frac{\Gamma(c)\Gamma(c-a-b_1)}{\Gamma(c-a)\Gamma(c-b_1)}\,\,_{2}\mathrm{F}_{1}\left( \left. 
\begin{array}{c}
a;b_2 \\[2mm]
c-b_1
\end{array}
\right| x\right)
\]
In similar way, it is possible to recognize the following formula, provided that $c>a+b_3$, which we used to simplify equation \eqref{solivp1k}.
\begin{equation}\label{toeta0b}
\mathrm{F}_{D}^{(3)}\left( \left. 
\begin{array}{c}
a;b_1,b_2,b_3 \\[2mm]
c
\end{array}
\right|x,y,1\right)=\frac{\Gamma(c)\Gamma(c-a-b_3)}{\Gamma(c-a)\Gamma(c-b_3)}\,\mathrm{F}_{1}\left( \left. 
\begin{array}{c}
a;b_1,b_2 \\[2mm]
c-b_3
\end{array}
\right| x,y\right)
\end{equation}
We also recall the reduction formula\footnote{ http://functions.wolfram.com/HypergeometricFunctions/AppellF1/03/05/}:
\begin{equation}\label{wolfram}
\mathrm{F}_{1}\left( \left. 
\begin{array}{c}
a;b,b, \\[2mm]
c
\end{array}
\right|x,-x\right)=\,_3\mathrm{F}_{2}\left( \left. 
\begin{array}{c}
\frac{a+1}{2},\frac{a}{2},b \\[2mm]
\frac{c+1}{2},\frac{c}{2}
\end{array}
\right|x^2\right)
\end{equation}

\subsubsection*{Acknowledgements}
The authors are indebted to their friend Prof. Aldo Scimone who drew some
figures of this paper: they hereby take the opportunity to thank him warmly. The second author is supported by RFO Italian grant funding.


\begin{thebibliography}{10}

\bibitem{Basoco}
M.~Basoco.
\newblock On the inflexional elastica.
\newblock {\em American Mathematical Monthly}, pages 303--309, 1941.

\bibitem{bele}
V.~Beletsky and E.~Levin.
\newblock {\em Dynamics of space tether systems}, volume~83.
\newblock Univelt Incorporated, San Diego, 1993.

\bibitem{Bellu}
O.~Belluzzi.
\newblock {\em Scienza delle costruzioni}, volume~I.
\newblock Zanichelli, Bologna, 1970.

\bibitem{burgatti}
P.~Burgatti.
\newblock {\em Teoria matematica dell'elasticit\`{a}}.
\newblock Zanichelli, Bologna, 1931.

\bibitem{turzi}
R.~De~Pascalis, G.~Napoli, and S.~Turzi.
\newblock Growth-induced blisters in a circular tube.
\newblock {\em Physica D}, 283:1--9, 2014.

\bibitem{dutka}
J.~Dutka.
\newblock The early history of the hypergeometric function.
\newblock {\em Archive for History of Exact Sciences}, 31(1):15--34, 1984.

\bibitem{1}
R.~Frisch-Fay.
\newblock {\em Flexible bars}.
\newblock Butterworths, London, 1962.

\bibitem{Goss}
V.~Goss.
\newblock The history of the planar elastica: insights into mechanics and
  scientific method.
\newblock {\em Science \& Education}, 18(8):1057--1082, 2009.

\bibitem{green}
G.~Greenhill.
\newblock Graphical representation of the elliptic functions by means of a bent
  elastic beam.
\newblock {\em Messenger of mathematics}, 5:180--188, 1876.

\bibitem{l}
G.~Lauricella.
\newblock Sulle funzioni ipergeometriche a pi\`{u} variabili.
\newblock {\em Rendiconti del Circolo Matematico di Palermo}, 7:111--158, 1893.

\bibitem{R1}
J.~A.~L. Matheson.
\newblock {\em Hyperstatic structures: an introduction to the theory of
  statically indeterminate structures}.
\newblock Butterworths, London, 1971.

\bibitem{el1}
G.~Mingari~Scarpello and D.~Ritelli.
\newblock Elliptic integral solutions of spatial elastica of a thin straight
  rod bent under concentrated terminal forces.
\newblock {\em Meccanica}, 41(5):519--527, 2006.

\bibitem{el2}
G.~Mingari~Scarpello and D.~Ritelli.
\newblock Elliptic integrals solution to elastica's boundary value problem of a
  rod bent by axial compression.
\newblock {\em Journal of Analysis and Applications}, 5(1):53--69, 2007.

\bibitem{el3}
G.~Mingari~Scarpello and D.~Ritelli.
\newblock Elliptic functions solution to exact curvature elastica of a thin
  cantilever under terminal loads.
\newblock {\em Journal of Geometry and Symmetry in Physics}, 12(1):75--92,
  2008.

\bibitem{Hindawi}
G.~Mingari~Scarpello and D.~Ritelli.
\newblock Exact solutions of nonlinear equation of rod deflections involving
  the Lauricella hypergeometric functions.
\newblock {\em International Journal of Mathematics and Mathematical Sciences},
  2011, 2011.

\bibitem{s}
S.~Saran.
\newblock Hypergeometric functions of three variables.
\newblock {\em Ganita}, 5:77--91, 1954.

\bibitem{Schell}
W.~Schell.
\newblock {\em Theorie der Bewegung und der Kr{\"a}fte}.
\newblock B.G. Teubner, Leipzig, 1870.

\bibitem{Tric}
F.~G. Tricomi.
\newblock {\em Funzioni ellittiche}.
\newblock Zanichelli, Bologna, 1951.

\bibitem{wag}
T.~Wagner and D.~Vella.
\newblock The sticky elastica: the lamination blisters beyond small
  deformations.
\newblock {\em Soft Matter}, 9:1025--130, 2013.

\end{thebibliography}

\end{document}